%% This is an Amstex file. The first version was
%% of April 19, 22, 24, 25, 30, 2007. May 2, 5, 2007.
%% Revised in Fisciano, Italy, June 4, 5, 2007. Again July 8, 2007.
%% The second version Bugeaud-Ivic of July 14, 15, 2007, was
%% file BugIvic1.tex, first version was BugdIvic.tex.
%% This is the version BugIvic4.tex of July 24-26, 2007,
%% after BugIvic2.tex, BugIvic3.tex of July 18-23, 2007.
%% Not for publication or circulation!

\input amstex.tex
\documentstyle{amsppt}
\magnification1200
\hsize=12.5cm
\vsize=18cm
\hoffset=1cm
\voffset=2cm

\def\DJ{\leavevmode\setbox0=\hbox{D}\kern0pt\rlap
{\kern.04em\raise.188\ht0\hbox{-}}D}
\def\dj{\leavevmode
 \setbox0=\hbox{d}\kern0pt\rlap{\kern.215em\raise.46\ht0\hbox{-}}d}

\def\txt#1{{\textstyle{#1}}}
\baselineskip=13pt
\def\hf{{\textstyle{1\over2}}}

\def\d{{\,\roman d}}
\def\e{\varepsilon}

\def\G{\Gamma}

\def\={\;=\;}

\def\zt{\zeta(\hf+it)}

\def\D{\Delta}

\def\R{\Re{\roman e}\,} 
\def\z{\zeta}

\def\hf{{\textstyle{1\over2}}}
\def\txt#1{{\textstyle{#1}}}

%%%%%%%%%%% Fonts macros %%%%%%%%%%%%
\font\tenmsb=msbm10
\font\sevenmsb=msbm7
\font\fivemsb=msbm5
\newfam\msbfam
\textfont\msbfam=\tenmsb
\scriptfont\msbfam=\sevenmsb
\scriptscriptfont\msbfam=\fivemsb
\def\Bbb#1{{\fam\msbfam #1}}

\def \NN {\Bbb N}

\def \ZZ {\Bbb Z}

\font\ff=cmr8
\def\txt#1{{\textstyle{#1}}}
\baselineskip=13pt

\font\teneufm=eufm10
\font\seveneufm=eufm7
\font\fiveeufm=eufm5
\newfam\eufmfam
\textfont\eufmfam=\teneufm
\scriptfont\eufmfam=\seveneufm
\scriptscriptfont\eufmfam=\fiveeufm
\def\mathfrak#1{{\fam\eufmfam\relax#1}}

\font\tenmsb=msbm10
\font\sevenmsb=msbm7
\font\fivemsb=msbm5
\newfam\msbfam
     \textfont\msbfam=\tenmsb
      \scriptfont\msbfam=\sevenmsb
      \scriptscriptfont\msbfam=\fivemsb
\def\Bbb#1{{\fam\msbfam #1}}

\def \NN {\Bbb N}

\def \ZZ {\Bbb Z}

  \def\rightheadline{{\hfil{\ff
Sums of the error term function in the mean square for
$\z(s)$}\hfil\tenrm\folio}}

  \def\leftheadline{{\tenrm\folio\hfil{\ff
  Yann Bugeaud and Aleksandar Ivi\'c }\hfil}}
  \def\emptyheadline{\hfil}
  \headline{\ifnum\pageno=1 \emptyheadline\else
  \ifodd\pageno \rightheadline \else \leftheadline\fi\fi}

\font\ff=cmr8
\font\teneufm=eufm10
\font\seveneufm=eufm7
\font\fiveeufm=eufm5
\newfam\eufmfam
\textfont\eufmfam=\teneufm
\scriptfont\eufmfam=\seveneufm
\scriptscriptfont\eufmfam=\fiveeufm
\def\mathfrak#1{{\fam\eufmfam\relax#1}}

\font\tenmsb=msbm10
\font\sevenmsb=msbm7
\font\fivemsb=msbm5
\newfam\msbfam
\textfont\msbfam=\tenmsb
\scriptfont\msbfam=\sevenmsb
\scriptscriptfont\msbfam=\fivemsb
\def\Bbb#1{{\fam\msbfam #1}}

\def \NN {\Bbb N}

\def \ZZ {\Bbb Z}

\def\D{\Delta}

 \def\e{\varepsilon}
 \def\d{\,{\roman d}}
%\topglue1cm
\topmatter
\title
Sums of the error term function in the mean square for $\z(s)$
\endtitle
\author
Yann Bugeaud and Aleksandar Ivi\'c
\endauthor
\address
Universit\'e Louis Pasteur, Math\'ematiques, 7 rue Ren\'e Descartes,
F-67084 Strasbourg cedex, France
\endaddress
\address
Katedra Matematike RGF-a,
Universitet u Beogradu,  \DJ u\v sina 7,
11000 Beograd, Serbia.\bigskip
\endaddress
\subjclass 11 M 06
\endsubjclass
\email
{\tt bugeaud\@math.u-strasbourg.fr,
ivic\@rgf.bg.ac.yu}
\endemail
\dedicatory
\enddedicatory
\abstract Sums of the form $\sum_{n\le x}E^k(n)\;(k\in\NN$ fixed)
are investigated, where $$
E(T) = \int_0^T|\zt|^2\d t - T\Bigl(\log {T\over2\pi} + 2\gamma -1\Bigr)
$$
is the error term in the mean square formula for $|\zt|$. The emphasis
is on the case $k=1$, which is more difficult than the corresponding
sum for the divisor problem. The analysis requires  bounds for the
irrationality measure of
${\roman e}^{2\pi m}$ and  for the partial quotients
in its continued fraction expansion.
\endabstract
\endtopmatter
\document
\head
1. Introduction
\endhead
Recently J. Furuya [3] investigated
 the sums of $\Delta^k(n)\;(n\in\NN)$ when $k=2$ or
$k=3$, and as usual
$$
\D(x) = \sum_{n\le x} d(n) - x(\log x + 2\gamma-1)
$$
denotes the error term in the classical Dirichlet divisor problem
($d(n)$ is the number of divisors of $n$, and $\gamma = -\G'(1) = 0.577215649\ldots$
is Euler's constant). For an account of the divisor problem, the reader
is referred to [7] and [8].
Furuya proved that
($c_1,c_2$ are suitable constants)
$$
\sum_{n\le x}\Delta^2(n) = \int_1^x \Delta^2(u)\,{\d}u + {1\over6}x\log^2x
+ c_1x\log x +  c_2x + R(x),
$$
where $R(x) = O(x^{3/4}\log x)$ and at the same time
$R(x) =\Omega_\pm(x^{3/4}\log x)$.
Here and below, $f(x) = \Omega_\pm(g(x))$ means that both
$$\limsup_{x\to\infty }f(x)/g(x)>0, \quad
\liminf_{x\to\infty} f(x)/g(x)<0$$
hold. This improves much on G.H. Hardy's
classical result [5] that $R(x) = O_\e(x^{1+\epsilon})$, where the
subscript means that the $O$--constant depends only on $\e$, an arbitrarily
small positive constant. The work [3] also brings forth a corresponding
asymptotic formula for $\sum_{n\le x}\Delta^3(n) $, where again the
error term is precisely determined. It seems natural to try to
generalize Furuya'a results to sums of $E^k(n)$, where $k\in\NN$ is fixed,
$\z(s) = \sum_{n=1}^\infty n^{-s}\;(\R s > 1)$ is the Riemann zeta-function, and
$$
E(T) = \int_0^T|\zt|^2\d t - T\Bigl(\log {T\over2\pi} + 2\gamma -1\Bigr)
\leqno(1.1)
$$
is the error term in the mean square formula for $|\zt|$. This is motivated
by the general analogy between $\D(x)$ and $E(T)$, inspired by the pioneering
work of F.V. Atkinson [1], who produced an explicit formula for $E(T)$ with
a good error term (see also [7] and [8]). Later research established more
precise analogies between the two functions, and in particular we refer to the
works of M. Jutila [14], [15], and the second author [10], [11].

\medskip In considering the sums $\sum_{n\le x}E^k(n)$, it turns out
that already for $k=1$ the problem is somewhat different than the
corresponding problem for sums of $\D(n)$.
Namely it reduces to the evaluation of
$$
\int_0^T\psi(t)|\zt|^2\d t\qquad\Bigl(\psi(t) = t - [t] - \hf\Bigr).\leqno(1.2)
$$
To see this, note that
$$
\eqalign{
\sum_{n\le x} E(n)&= \sum_{n\le x} (E(n)-\pi) + \pi x +O(1)\cr&
= \pi x +O(1) + \int_0^x(E(u)-\pi)\d[u]\cr&
= \pi x +O(1) + G(x)  - \int_0^x(E(u)-\pi)\d\psi(u)\cr&
= \pi x +O(\log x) + G(x) + O(|E(x)|) + \int_0^x\psi(t)|\zt|^2\d t,\cr}
$$
where
$$
G(x) := \int_0^x(E(u)-\pi)\d u,\quad G(x) = O(x^{3/4}), \quad
G(x) = \Omega_\pm(x^{3/4}).\leqno(1.3)
$$
The results on $G(x)$ in (1.3) are to be found in the second author's joint
work with  J.L. Hafner [4].
In what concerns the summatory function of $\D(n)$,
we recall that already G.F. Vorono{\"\i} [19] proved that
$$
\sum_{n\le x}\Delta(n) = \hf x\log x + (\gamma - {\txt{1\over4}})x + R_0(x),\leqno(1.4)
$$
where $R_0(x) = O(x^{3/4})$ and at the same time $R_0(x)= \Omega_\pm(x^{3/4})$.

\medskip
 In working out the proof of (1.4), one actually obtains
$$
\sum_{n\le x}\Delta(n) = \hf x\log x + (\gamma - {\txt{1\over2}})x + \D(x)
+ \int_1^x\D(t)\d t + O(\log x),\leqno(1.5)
$$
and (1.4) follows from (1.5) when one recalls Vorono{\"\i}'s classical results [19] that
$$
\int_1^x \D(t)\d t = {x\over4} + R_1(x),\quad R_1(x) = O(x^{3/4}), \quad
R_1(x) = \Omega_\pm(x^{3/4}).\leqno(1.6)
$$
The results in (1.3) on $G(x)$ are in fact the analogues of Vorono{\"\i}'s
results in (1.6).

\medskip
Let us point out the
basic difference between the (discrete) function $\D(n)$ and the (continuous)
function $E(n)$ in the context of their respective summatory
functions. For the former we start (we may clearly assume that $x\in \NN$
in what follows) with
$$
\sum_{n\le x}\Delta(n) = \sum_{m\le x}d(m)
\Bigl(\sum_{m\le n\le x}1\Bigr) - \sum_{n\le x}n(\log n + 2\gamma-1). \leqno(1.7)
$$
The last sum is easily evaluated by the Euler-Maclaurin summation formula. Since we
have
$$
\sum_{m\le n\le x}1 = x-m+1\qquad(m,\,x \in \NN),\leqno(1.8)
$$
one arrives without any major difficulties at (1.5). But in the analogue of (1.7) for
$E(n)$, the first sum on the right-hand side of (1.7) will correspond to
$$
\int_0^x|\zt|^2\Bigl(\sum_{t\le n\le x}1\Bigr)\d t.\leqno(1.9)
$$
Observe that in (1.9) $t$ is not an integer (except for $O(x)$
exceptions), and an identity like (1.8) will not hold for the sum
over $n$ in (1.9). Instead, the function $[t] = t  - \psi(t) - \hf$
will come into play, producing (1.2). It is rather curious, as will
be shown in the next section, that our problem will reduce to bounds
for the irrationality measure of ${\roman e}^{2\pi m}$ and  for the
partial quotients in its continued fraction expansion.

\medskip {\bf Acknowledgements}.
We wish to thank Prof. Michel Waldschmidt
for valuable remarks.

\bigskip
\head
2. The summatory function of $E(n)$
\endhead
In this section we shall continue the discussion on the summatory function of $E(n)$.
In Section 1 we established that
$$
\sum_{n\le x}E(n) = \pi x + \int_0^x\psi(t)|\zt|^2\d t + G(x) + O(x^{1/3}),
\leqno(2.1)
$$
where $G(x)$ is given by (1.3), and where we used the (weak) bound (see [7]
or [8]) $E(x) \ll x^{1/3}$. Since $\psi(t)$ is an oscillating function, one
expects that there will be a lot of cancellations in the integral in (2.1),
but this turns out to be hard to prove.
\medskip
To bound the integral (2.1), we first consider subintegrals
over intervals of the form $[T,\,2T]\,$. Then we use the
approximate functional equation for $\z^2(s)$ (see e.g., [7] and [10])
in the form
$$
|\zt|^2 = 2\sum_{n\le t/(2\pi)}d(n)n^{-1/2}
\cos\left(t\log{t\over2\pi n}-t-{\pi\over4}\right) +{\Cal R}(t),\leqno(2.2)
$$
where the error term function ${\Cal R}(t)$
satisfies the mean square relation
$$
\int_1^T{\Cal R}^2(t)\d t  \;=\; AT^{1/2} + O(\log^5T)\qquad( A>0).\leqno(2.3)
$$
The asymptotic formula (2.3) was proved by I. Kiuchi [16], who used in his proof precise
asymptotic expressions for ${\Cal R}(t)$ (related to $\D(t/(2\pi))$,
due to Y. Motohashi [17].
Using (2.3) we obtain, by the Cauchy-Schwarz inequality for integrals,
$$\eqalign{&
\int_T^{2T}\psi(t)|\zt|^2\d t \,= \,O(T^{3/4}) \;+ \cr&\,+ \sum_{n\le T/\pi}d(n)n^{-1/2}
\int_{\max(T,2\pi n)}^{2T}\psi(t)
\cos\left(t\log{t\over2\pi n}-t-{\pi\over4}\right)\d t.
\cr}
$$
For the integral on the right-hand side  we use the
familiar Fourier series
$$
\psi(x) = - {1\over\pi}\sum_{m=1}^\infty {\sin(2\pi mx)\over m}\qquad(x\not\in \ZZ),
\leqno(2.4)
$$
which is boundedly convergent and thus may be integrated termwise.
This leads to
$$\eqalign{
&\int_T^{2T}\psi(t)|\zt|^2\d t\; = \; O(T^{3/4})\;-\cr&
-{1\over\pi}\sum_{m\le\log T}{1\over m}
\sum_{n\le T/\pi}d(n)n^{-1/2}
\int\limits_{\max(T,2\pi n)}^{2T}\sin(2\pi mt)
\cos\left(t\log{t\over2\pi n}-t-{\pi\over4}\right)\d t.\cr}
$$
Namely, when we write the sines and cosines as exponentials,
the relevant exponential integrals will be of the form
$$
\int_{\max(T,2\pi n)}^{2T}{\roman e}^{\pm iF(t)}\d t\quad
\qquad\left(1\le n\le {T\over\pi}\right),
$$
where
$$
F(t) =t\log{t\over2\pi n}-t-2\pi mt,
\; F'(t) = \log{t\over2\pi n} - 2\pi m,\; F''(t) = {1\over t}.
$$
Therefore for $m>\log T$ we have $|F'(t)|\gg m$, hence by the
first derivative test for exponential integrals (see e.g., Lemma 2.1
of [7]) the contribution of such $m$ is $O(T^{1/2}\log T)$. Next,
for $T/(2\pi) \le n \le T/\pi, $ we have
$$
|F'(t)| = 2\pi m - \log{t\over2\pi n} \ge 2\pi m - \log {t\over T}
\ge 2\pi m - \log 2 \ge \pi m\quad(m\in\NN),
$$
hence the contribution of such $n$ is again $O(T^{1/2}\log T)$.
Likewise the contribution of satisfying $m > (1/(2\pi)+\e)\log T$ is
$O(T^{3/4})$. What remains is
$$
\eqalign{
&\int_T^{2T}\psi(t)|\zt|^2\d t \;=\;  O(T^{3/4})\cr
-{1\over\pi}\sum_{m\le(1/(2\pi)+\e)\log T}{1\over m} &
\sum_{n\le T/(2\pi)}{d(n)\over n^{1/2}}
\int\limits_{T}^{2T}\sin(2\pi mt)
\cos\left(t\log{t\over2\pi n}-t-{\pi\over4}\right)\d t.
\cr} \leqno(2.5)
$$

The saddle point (the root of $F'(t) = 0$)
is $t_0 = 2\pi n{\roman e}^{2\pi m}$, and it  lies in $[T,\,2T]\,$ for
$n \asymp T{\roman e}^{-2\pi m}$. Using the saddle-point method
(direct integration over a suitable contour in the complex plane
joining the points $T$, $2T$ and passing through $t_0$)
we obtain   that the expression in (2.5) is asymptotic to
sums of the type ($C>0$ is a generic constant)
$$
C\sum_{m\le (1/(2\pi)+\e)\log T}{{\roman e}^{\pi m}\over m}
\sum_{n\asymp T{\roman e}^{-2\pi m}}d(n)
\exp\left(-2\pi in{\roman e}^{2\pi m}\right).\leqno(2.6)
$$
%% ai
%% Pas d'erreur, on prend la valeur conjugee, c'est la meme
%% chose parce qu'on veut estimer la valeur absolue.

The problem is thus reduced to the estimation of the sum over $n$ in (2.6).
One can use the results of a classic paper by J.R. Wilton [21].
He proved, among other things, the functional equation
$$
D(x,\eta) = \eta^{-1}D(\eta^2x,-\eta^{-1}) + O(x^{1/2}\log x),
\quad D(x,\eta) = \sum_{n\le x}d(n){\roman e}^{2\pi i\eta n},
\leqno(2.7)
$$
where $\eta^2 x\gg 1$ and $0< \eta\le 1$ is real.
After taking the conjugate sum in (2.6), (2.7) can be applied with
 $x = CT{\roman e}^{-2\pi m},\, \eta = \{{\roman e}^{2\pi m}\}$,
 where $\{y\} = y - [y]$ denotes the fractional part of $y$.
The problem is that it is difficult to find an explicit expression
for $\eta = \eta(m)$. On the other hand, Wilton (op. cit.)
states that for irrational $\eta$ one always has
$$
D(x,\eta) \;=\;o(x\log x)\qquad(x\to\infty)\leqno(2.8)
$$
uniformly in $\eta = \eta(m)$. The bound (2.8) can be applied in our case,
since the numbers $\eta(m)$ are all transcendental.
Namely, since the number ${\roman e}^\pi$ is
transcendental,  all the numbers ${\roman e}^{2\pi m}$
are also transcendental.
The fact that ${\roman e}^\pi = (-1)^{-i}$
is transcendental follows from the solution of Hilbert's 7th problem
(see e.g., C.L. Siegel [18] for a proof) that $\alpha^\beta$
is transcendental for algebraic $\alpha \neq 0,1$ and
irrational $\beta$. This was proved by Gel'fond and Schneider
independently in 1934. Wilton [21] also proved
the following: let
$$
{1\over\eta} = a_0 + {1\over a_1+}{1\over a_2+ }\cdots\qquad(0<\eta \le 1)
$$
be the expansion of $1/\eta$ as a simple continued fraction.
If $\eta$ is irrational and the convergents $a_n$ satisfy
$a_n \ll n^{1+K}\;(K\ge0)$, then
$$
\sum_{n\le x}d(n)\sin(2\pi n\eta) \ll x^{1/2}\log^{2+K}x,\leqno(2.9)
$$
while if
$$
a_n \ll {\roman e}^{Kn}\;(K>0)\leqno(2.10)
$$
holds, then
$$
\sum_{n\le x}d(n)\exp(2\pi in\eta) \ll x^{H}\log x
\quad\left(H = {4K+\log2\over 4K+ 2\log 2}\right).\leqno(2.11)
$$
Therefore the information on the size of the $a_n$'s would lead, by
(2.9) or (2.11), to good bounds for $\sum_{n\le x}d(n)\sin(2\pi n\eta)$,
which is the sum we need to estimate. Additional care should be exerted
because in our case $\eta = {\roman e}^{2\pi m}$, so that $a_n = a_n(m)$
will depend on $m$ as well. To determine the true order of $a_n = a_n(m)$
seems to be a difficult problem; even the fact that ${\roman e}^{\pi}$ is not
a Liouville number is not known. In the next section we shall provide
an explicit bound for the partial quotients of ${\roman e}^{2\pi m}$.
Although poorer than (2.10),
it is still sufficient to provide a non-trivial bound which in our case improves
on the general bound (2.8) in the case when $\eta$ is irrational.

\head
3. The results on the sum of $E(n)$
\endhead
\smallskip

We begin by recalling a bound for  the irrationality measure
of ${\roman e}^{\pi m}$. Namely,
for any positive integers $m$, $p$ and $q$ with $p \ge 3$, we have
$$
\biggl| {\roman e}^{\pi m} - {p \over q} \biggr| >
\exp\Bigl\{-2^{72} \,  (\log 2m)
(\log p) (\log\log p)\Bigr\}.\leqno(3.1)
$$
This result was proved by M. Waldschmidt [20, p. 473].
Now for any non-zero
integer $m$ let
$$
{\roman e}^{\pi m} = [a_0(m) ; a_1(m), \ldots ]
$$
be the expansion of ${\roman e}^{\pi m}$ as a continued fraction.

\bigskip
LEMMA 1.
{\it There exists an absolute positive constant $c$ such that
$$
\log\log a_n (m)
< c (n+\log m) \log (n + \log m) \leqno(3.2)
$$
holds for any non-zero integer $m$ and any non-negative integer $n$.}

\bigskip

It is sufficient to establish Lemma 1 for positive values of $m$,
since ${\roman e}^{-\pi m} = [0; a_0(m), a_1(m), \ldots ]$
for $m \ge 1$.
Let $\xi$ be a positive real number with continued fraction expansion
$$
\xi \;= \;[a_0; a_1, a_2, \ldots]
$$
and set, for $n \ge 1$,
$$
{p_n \over q_n} \;=\; [a_0; a_1, a_2, \ldots , a_n].
$$
From the theory of continued fractions it is known
that
$$
\biggl| \xi - {p_n \over q_n} \biggr|
< {1 \over q_n q_{n+1}} \qquad (n \ge 1).  \leqno (3.3)
$$
It follows from (3.1) that
$$
\biggl| {\roman e}^{\pi m} - {p \over q} \biggr| >
\exp\{-2^{75} (\log 2m) (\log q)(\log\log q)\},
\leqno (3.4)
$$
for all $p/q$ with $q\ge {\roman e}^{\pi m}$,
since $p \le 2 q {\roman e}^{\pi m}$, giving
$$
\log p \le \log 2 + \log q + \pi m \le \log 2 + 2\log q
\le 3 \log q.
$$
In particular, (3.4) is true for the convergents
$p_h(m)/ q_h(m)$ to ${\roman e}^{\pi m}$ as soon as
$q_h (m) \ge {\roman e}^{\pi m}$. Let $n \ge 2$
be an integer. Writing simply $p_n/q_n$
instead of $p_n(m) / q_n (m)$, it then follows
from (3.3) and (3.4) that
$$
q_n <  \exp\Bigl\{2^{75} (\log 2m) (\log q_{n-1})
(\log\log q_{n-1})\Bigr\},  \leqno (3.5)
$$
if $q_{n-1} \ge {\roman e}^{\pi m}$.
We infer from (3.5) that there exists an absolute
positive constant $c_1$ such that
$$
\log_2 q_n < c_1 + \log_2 m + \log_2 q_{n-1}
+ \log_3 q_{n-1}, \leqno (3.6)
$$
where, for $j\ge2$, $\,\log_jx = \log(\log_{j-1}x)\;(\log_1x \equiv \log x)$
denotes the $j$-th iteration of the logarithm.
On dividing (3.6) by $\log_3 q_n$, we obtain
$$
{\log_2 q_n \over \log_3 q_n} < 1 +
(c_1 + \log_2 m) \, {1 \over \log_3 q_n}  +
{\log_2 q_{n-1} \over \log_3 q_n}.  \leqno (3.7)
$$
Using repeatedly the fact that
$q_n > q_{n-1}$, it follows that
$$
\eqalign{
{\log_2 q_n \over \log_3 q_n} & < n +
(c_1 + \log_2 m) \, \sum_{j=h+1}^n \, {1 \over \log_3 q_j}
+ {\log_2 q_{h} \over \log_3 q_{h+1}} \cr
& <   n +
n \, {c_1 + \log_2 m \over \log_3 q_h} \, +
\log_2 q_h, \cr}  \leqno (3.8)
$$
where $h$ is the integer such that $q_{h-1} < e^{\pi m}$
and $q_{h} > e^{\pi m}$.
This choice of $h$ and (3.1) imply that
$$
q_{h} < \exp\{c_2 m (\log 2m)^2\},  \leqno (3.9)
$$
where $c_2$ (as $c_3$, $c_4$ and $c_5$ below) is
an absolute positive constant.
Consequently, we infer from (3.8) and (3.9) that
$$
{\log_2 q_n \over \log_3 q_n} < n +
n \, {c_1 + \log_2 m \over \log_2 m} +
c_3 + \log m
< c_4 (n + \log m).
$$
This gives at once the asserted upper bound (3.2), namely
$$
\log_2 a_n < \log_2 q_n
< c_5 (n+\log m) \log (n + \log m).
$$
Now we rewrite (3.2) as
$$
a_n(m) < \exp\left\{(n+\log m)^{c(n+\log m)}\right\}\leqno(3.10)
$$
and proceed to estimate, in view of the transformation formula (2.7),
$$
D(x,\eta(m)) := \sum_{n\le x}d(n)\exp(2\pi i\eta(m))\qquad\Bigl(\eta(m)
= {\roman e}^{-2\pi m},
\;1 \le m \le {\log_2x\over\log_3x}\Bigr).\leqno(3.11)
$$
Then $1/\eta(m) = {\roman e}^{2\pi m}$, and the above discussion applies with $m$
replaced by $2m$.
The above range for $m$ suffices, since by trivial estimation the contribution of
$m > \log_2T/\log_3T$ to (2.6) is
$$
\ll \sum_{m>\log_2T/\log_3T}{{\roman e}^{\pi m}\over m}\cdot {T\over {\roman e}^{2\pi m}}
\log T \ll T\log T\exp\Bigl\{-C{\log_2T\over\log_3T}\Bigr\}.
$$
From Wilton's bounds [21, eqs. $(6.3_{III})$ and $(6.3_{IV})$] we have
$$
D(x,\eta(m)) \ll x^{1/2}\log^2x + \min\Bigl(a_N(2m)x^{1/2}\log x,\,
2^{-{1\over2}N}x\log x\Bigr),\leqno(3.12)
$$
for large $N$ satisfying $N \ll \log x$.
We use (3.10) to deduce that, for
$$
N < \Bigl({1\over c}+o(1)\Bigr){\log_2x\over\log_3x},
$$
 we have
$$\eqalign{
a_N(2m)x^{1/2}\log x &\ll x^{1/2}\log x\exp\left\{\Bigl(N+{\log_2 x\over\log_3x}\Bigr)^
{c(N+\log_2 x/\log_3x)}\right\}
\cr&
\ll x\log x\exp\left(-C{\log_2x\over\log_3x}\right)\quad(C = C(c)>0),
\cr}
$$
while for $N \ge (1/c+o(1))\log_2x/\log_3x$ obviously
$$
2^{-{1\over2}N}x\log x \ll x\log x\exp\left(-C{\log_2x\over\log_3x}\right).
$$
This means that we obtain the following non-trivial result, which improves on (2.8),
although only by a quantity that is less than $(\log x)^\e$ for any given $\e>0$.
\bigskip
THEOREM 1. {\it If $D(x,\eta(m))$ is given by} (3.11), {\it then for
$1\le m \le {\log_2x\over\log_3x}$ and some constant $C>0$ we have}
$$
D(x,\eta(m))\;\ll\;x\log x\exp\left(-C{\log_2x\over\log_3x}\right).\leqno(3.13)
$$

\bigskip

To apply Theorem 1 we transform the sum over $n$ in (2.6) by (2.7) (in the
range $1\le m \le \log_2T/\log_3T$) with
$x = CT{\roman e}^{-2\pi m}$, $\eta = {\roman e}^{2\pi m}$. We obtain
$D(x,\eta(m))$ with $\eta(m) = {\roman e}^{-2\pi m}, x = CT{\roman e}^{2\pi m}$,
plus an error term of $O(T^{1/2}\log T)$.
We use (3.13)  for this $D(x,\eta(m))$, estimate the contribution of the
range $ m \ge \log_2T/\log_3T$ trivially, and arrive at

\bigskip
THEOREM 2. {\it We have
$$
\sum_{n\le x} E(n) \;=\; \pi x + H(x),\leqno(3.14)
$$
where, for some $C>0$, unconditionally
$$
H(x) \;\ll\;x\log x\exp\left(-C{\log_2x\over\log_3x}\right).\leqno(3.15)
$$
If} (2.9) {\it holds with $\eta = {\roman e}^{2\pi m}$ uniformly
for $1\le m\ll \log x$, then we have}
$$
H(x) \;=\;O(x^{3/4}).\leqno(3.16)
$$

\medskip
{\bf Remark 1}. The intrinsic difficulty in this problem is to bound
the sum $D(x,\eta(m))$.
In fact, unconditionally our bound for $H(x)$ in (3.15) is poorer
than the main term $\pi x$ in (3.14). However,
 if the
bound in (2.9) is true in our case, then (3.16) can be complemented
with $H(x) = \Omega_\pm(x^{3/4})$. Namely instead of using
(2.2)--(2.3), we may use an (unsymmetric) approximate functional
equation for $|\zt|^2$ with a sharp error term (see e.g., [7,
Chapter 4] ). In fact it may be conjectured that instead of (3.15)
with (3.16) one has
$$
\sum_{n\le x} E(n) \;=\; \pi x + G(x) + O_\e(x^{1/2+\e}),
$$
where $G(x)$ is given by (1.3).

\medskip
{\bf Remark 2}. Another discrete sum involving the divisor function
was investigated by Coppola--Salerno [2] and A. Ivi\'c [12].
The former have shown that, for $T^\e \le h \le \hf \sqrt{T},\,L = \log T$,
$$
\sum_{T\le n \le 2T}\left(\D(n+U)-\D(n)\right)^2=
{8\over\pi^2}TU\log^3\Bigl({\sqrt{T}\over U}\Bigr)
+ O(TUL^{5/2}\sqrt{L}).\leqno(3.17)
$$
This is an asymptotic formula with a weak error term, improved in [12]
by the second author to
$$\eqalign{
\sum_{T\le n \le 2T}\left(\D(n+U)-\D(n)\right)^2 & = TU\sum_{j=0}^3c_j\log^j
\Bigl({\sqrt{T}\over U}\Bigr) \cr&
+ O_\e(T^{1/2+\e}U^2) + O_\e(T^{1+\e}U^{1/2}),\cr}\leqno(3.18)
$$
where $1 \ll U \le \hf \sqrt{T}$ with $c_3 = 8\pi^{-2}$.
However, the analogue of
(3.17) or (3.18) for the sum
$$
\sum_{T\le n \le 2T}\left(E(n+U)-E(n)\right)^2\leqno(3.19)
$$
does not carry over, for the same reason for which we had difficulties in
evaluating sums of $E(n)$; because $E(T)$  is a continuous function,
while $\D(x)$ is not, having jumps at natural numbers of order at most
$O_\e(x^\e)$. The true order of magnitude of the sum in (3.19) seems elusive.

\head
4. Sums of $E^k(n)$ for $k>1$
\endhead
We begin by  noting that in general, for $k\in \NN$ fixed, one has
$$
\sum_{n\le x} E^k(n) = \int_0^x E^k(u)\d[u] = \int_0^x E^k(u)\d u -
\int_0^x E^k(u)\d\psi(u).\leqno(4.1)
$$
Integration by parts gives
$$
\eqalign{&\int_0^x E^k(u)\d\psi(u) = E^k(x)\psi(x) - k\int_0^x
\psi(u)E^{k-1}(u)E'(u)\d u\cr& \ll 1 + |E(x)|^k + \int_{10}^x
|E(u)|^{k-1}\left(|\z(\hf+iu)|^2 + \log\bigl({u\over2\pi}\bigr)+2\gamma\right)\d u,
\cr}
\leqno(4.2)
$$
where we used the defining property (1.1). Here we distinguish two cases.

\medskip
If $k = 2,4$, then by the Cauchy-Schwarz inequality for integrals the
last integral above is
$$
\eqalign{&
\ll \left\{\int_{10}^x E^{2k-2}(u)\d u\cdot
\int_{10}^x\Bigl(|\z(\hf+iu)|^4+\log^2u\Bigr)
\d u\right\}^{1/2}
\cr&
\ll x^{k/3} + \Bigl(x^{1+\hf(k-1)}\cdot x\log^4x\Bigr)^{1/2}
\ll x^{k/3} + x^{(k+3)/4}\log^2x
\cr& \ll x^{(k+3)/4}\log^2x,\cr}
$$
where we used (see the works of  D.R. Heath-Brown [6],
Ivi\'c--Sargos [13] and W. Zhai [22]),
$$
\int_0^x |E(u)|^k \d u \;\ll\; x^{1+k/4}\qquad(k\in \NN,\; 1\le k
\le 9).\leqno(4.3)
$$
Larger even values of $k$ could be, of course, also considered, but
for such values we do not have precise formulas for the integrals of $E^k$,
as we have when $k=2$ or $k=4$.
Namely we have (see [8])
$$
\int_0^x E^2(u)\d u = C_2x^{3/2} + O(x\log^4x)\qquad(C_2>0),\leqno(4.4)
$$
and (see W. Zhai [22, Part III])
$$
\int_0^x E^4(u)\d u = C_4x^{2} + O_\e(x^{2-3/28+\e})\qquad(C_4>0).\leqno(4.5)
$$

\medskip
If $k = 3,5,7,9$ then  $k-1$  is even, and we have
$|E(u)|^{k-1} = E^{k-1}(u)$. Since (4.3) holds and $E(x) \ll x^{1/3}$,
we see that the last integral in (4.2) is
$$
\int_{10}^xE^{k-1}(u)\Bigl(E'(u) + O(\log u)\bigr)\d u
\ll |E(x)|^k + x^{(k+3)/4}\log x \ll x^{k/3}\log x.
$$

\medskip
Note that we also have (see W. Zhai [23], who improved
the exponent $7/4 - 1/12$ of Ivi\'c--Sargos [13] in the error term)
$$
\int_0^x E^3(u)\d u = C_3x^{7/4} + O_\e(x^{7/4- 83/393+\e})\;\left(
{83\over393} = 0.211195\ldots\,,C_3>0\right),\leqno(4.6)
$$
and, in general (see W. Zhai [22]) for $5\le k\le9$,
$$
\int_0^x E^k(u)\d u = C_kx^{1+k/4} + O(x^{1+k/4-\delta(k)})\qquad(C_k>0),\leqno(4.7)
$$
where $\delta(k)$ is a positive
constant which may be explicitly evaluated.
From (4.1)--(4.7) it follows that we have proved

\bigskip
THEOREM 3. {\it We have
$$\eqalign{\cr&
\sum_{n\le x} E^2(n) = C_2x^{3/2} + O(x^{5/4}\log^2x),\cr&
\sum_{n\le x} E^3(n) = C_3x^{7/4} + O_\e(x^{7/4-83/393+\e}),\cr&
\sum_{n\le x} E^4(n) = C_4x^{2} + O_\e(x^{2-3/28+\e}),\cr}
$$
and for $k=5,7,9,$}
$$
\sum_{n\le x} E^k(n) = C_kx^{1+k/4} + O(x^{1+k/4-\delta(k)}) + O(x^{k/3}\log x).
$$

\medskip
Naturally, other values of $k$ could be also considered, but
as already remarked the results would
not be very good in view of the existing results on $\int_0^x E^k(u)\d u$.

\vfill
\eject
\topskip1cm

\Refs
\bigskip
\bigskip

\item{[1]} F.V. Atkinson, The mean value of the Riemann zeta-function,
Acta Math. {\bf81}(1949), 353-376.

\item{[2]} G. Coppola and S. Salerno, On the symmetry of the divisor
function in almost all short intervals, Acta Arith. {\bf113}(2004),
189-201.

\item{[3]} J. Furuya, On the average orders of the error
term in the Dirichlet divisor problems, J. Number Theory {\bf115}(2005), 1-26.

\item{[4]} J.L. Hafner and A. Ivi\'c,
On the mean square of the Riemann zeta-function on the critical line,
Journal of Number Theory {\bf 32}(1989), 151-191.

\item{[5]} G.H. Hardy, The average orders of the arithmetical functions
$P(x)$ and $\D(x)$, Proc. London Math. Soc. (2){\bf15}(1916), 192-213.

\item{[6]} D.R. Heath-Brown,
The distribution and moments of the error term in the Dirichlet divisor problems,
Acta Arith. {\bf60}(1992).

\item{[7]} A. Ivi\'c, The Riemann zeta-function, John Wiley \&
Sons, New York, 1985 (2nd ed. Dover, Mineola, New York, 2003).

\item{[8]} A. Ivi\'c, The mean values of the Riemann zeta-function,
LNs {\bf 82}, Tata Inst. of Fundamental Research, Bombay (distr. by
Springer Verlag, Berlin etc.), 1991.

\item{[9]} A. Ivi\'c,  Power moments of the error term in the
approximate functional equation for $\z^2(s)$, Acta Arith. {\bf65}(1993),
137-145.

\item{[10]} A. Ivi\'c, On the Riemann zeta-function and the divisor problem,
Central European J. Math. 2({\bf4}) (2004),   1-15;
II ibid. 3({\bf2}) (2005), 203-214.

\item{[11]} A. Ivi\'c, On the mean square of the zeta-function and
the divisor problem, Annales  Acad. Scien.
Fennicae Math. {\bf23}(2007), 1-9.

\item{[12]} A. Ivi\'c, On the divisor function and the Riemann
zeta-function in short intervals, to appear, {\tt arXiv:0707.1756.}

\item{[13]} A. Ivi\'c and P. Sargos,
On the higher moments of the error term in the divisor problem, Ill. J. Math. 2007,
in press, {\tt math.NT/0411537}, 27 pp.

\item{[14]} M. Jutila, Riemann's zeta-function and the divisor problem,
Arkiv Mat. {\bf21}(1983), 75-96 and II, ibid. {\bf31}(1993), 61-70.

\item{[15]} M. Jutila, On a formula of Atkinson, in ``Coll. Math. Sci.
J\'anos Bolyai 34, Topics in classical Number Theory, Budapest 1981",
North-Holland, Amsterdam, 1984, pp. 807-823.

\item{[16]} I. Kiuchi, An improvement on the mean value formula for the
approximate functional equation of the Riemann zeta-function,
J. Number Theory {\bf45}(1993), 312-319.

\item{[17]} Y. Motohashi,  A note on the approximate functional equation for
$\z^2(s)$, Proc. Japan Acad. Ser. A Math. Sci {\bf59}(1983), 393-396 and
II, ibid. {\bf59}(1983), 469-472.

\item{[18]} C.L. Siegel, Transcendental numbers, Annals of Mathematics studies,
no. {\bf16}, Princeton University Press, Princeton, N.J., 1949.

\item{[19]} G.F. Vorono{\"\i}, Sur une fonction transcendante et ses applications
\`a la sommation de quelques s\'eries, Ann. \'Ecole Normale (3){\bf21}(1904),
207-267 and ibid. 459-533.

\item{[20]} M. Waldschmidt, Simultaneous approximation of numbers connected
with the exponential function,
J. Austral. Math. Soc. Ser. A {\bf25}(1978), no. 4, 466-478.

\item{[21]} J.R. Wilton,
An approximate functional equation with applications to a problem of
Diophantine approximation, J. reine angew. Math. {\bf169}(1933), 219-237.

\item{[22]} W. Zhai,  On higher-power moments of $\D(x)$,
Acta Arith.  {\bf112}(2004),
367-395, II ibid. {\bf114}\break (2004), 35-54 and III ibid. {\bf118}(2005), 263-281.

\item{[23]} W. Zhai,  On higher-power moments of $E(t)$,
Acta Arith.  {\bf115}(2004), 329-348.

\vskip2cm
\endRefs

\enddocument

\bye